\documentclass[12 pt]{article}
\usepackage[english]{babel}
\usepackage{amsmath}
\usepackage{amsfonts}
\usepackage{amssymb}
\usepackage{cite}

\author{Philippe Jouan and Sa•d Naciri}
\title{Asymptotic Stability of Uniformly Bounded Nonlinear Switched Systems}

\newcommand\D{\, \textnormal{d}}
\newcommand\K{\mathcal{K}}
\newcommand\V{\mathcal{V}}
\newcommand\C{\mathcal{C}}
\newcommand\R{\mathbb{R}}

\newcommand\N{\mathbb{N}}
\newcommand\Li{\mathcal{L}}
\newcommand\M{\mathcal{M}}
\newtheorem{theo}{Theorem}
\newtheorem{cor}{Corollary}
\newtheorem{defi}{Definition}
\newtheorem{rema}{Remark}

\newtheorem{prop}{Proposition}
\newtheorem{lem}{Lemma}
\newtheorem{ex}{Example}

\numberwithin{equation}{section}
\newcommand{\demo}{\noindent\quad\textbf{\em Proof :$\,$} }
\newcommand{\findemo}{\begin{flushright}$\square$\end{flushright}} 

\begin{document}

\title{\bf Asymptotic Stability of Uniformly Bounded Nonlinear Switched Systems}
\author{Philippe JOUAN, Sa\"\i d NACIRI\footnote{LMRS, CNRS UMR 6085, Universit\'e
    de Rouen, avenue de l'universit\'e BP 12, 76801
    Saint-Etienne-du-Rouvray France. E-mail: Philippe.Jouan@univ-rouen.fr, Said.Naciri@etu.univ-rouen.fr}}

\date{\today}

\maketitle

\begin{abstract}
We study the asymptotic stability properties of nonlinear switched systems under the assumption of the existence of a common weak Lyapunov function.

We consider the class of nonchaotic inputs, which generalize the different notions of inputs with dwell-time, and the class of general ones. For each of them we provide some sufficient conditions for asymptotic stability in terms of the geometry of certain sets.

The results, which extend those of \cite{BaldeJouan}, are illustrated by many examples.    

\vskip 0.2cm


\vskip 0.2cm


\end{abstract}

\section{Introduction.}

This paper is concerned with the asymptotic stability properties of nonlinear switched systems defined by a finite collection $\{f_1,\ldots,f_p\}$ of smooth vector fields in $\R^d$. They are assumed to share a smooth weak Lyapunov function $V$, hence to be stable at the origin, and we are looking for sufficient conditions of asymptotic stability for certain classes of switching laws.\\
Our stability hypothesis is similar to the ones of \cite{Serre,BaldeJouan,BaldeJouan2,Bacciotti20051109,MancillaAguilar2006376}. 
In particular in \cite{BaldeJouan,BaldeJouan2,Serre} the vector fields are linear, the Lyapunov function is quadratic, and the asymptotic stability properties are closely related to the geometry of the union of some linear subspaces of $\R^d$.

In the present paper, we introduce in the same way two geometric subsets of $\R^d$ which turn out to be fundamental in the sense that they contain all the limit sets for two classes of inputs. More accurately let $\K_i=\{\Li_{f_i} V=0\}$ be the set where the vector field $f_i$ is tangent to a level set of the Lyapunov function $V$. We show in Section~\ref{moregeneral} that the limit sets are always contained in the union of the $\K_i$'s.  On the other hand, for the class of nonchaotic inputs (which are a generalization of inputs with dwell-time), the limit sets are included in the smaller set $\bigcup_{i=1}^p \M_i$ where $\M_i$ is roughly speaking the largest $f_i$-invariant subset of $\K_i$. 

The geometry of the union of the $\M_i$'s and the $\K_i$'s is consequently the key point of our results. 
Of course these sets are no longer linear whenever the vector fields are not, and we are lead to approximate them by linear subspaces or cones.

In the paper \cite{Bacciotti20051109} the authors deal with another set included in the union of the $\K_i$'s but generally much larger than the union of the $\M_i$'s (see Example~\ref{ridicule} in Section~\ref{Reguliere}).
However they consider only dwell-time inputs and their results are improved by the distinction we introduce between the two categories of sets and inputs.

The paper is organized as follows: basic definitions and notations are stated in Section~\ref{Position}, and some preliminary and technical results in Section~\ref{Convergence}.
In Section~\ref{Reguliere}, we state  a general result of asymptotic stability for nonchaotic inputs. 
We provide sufficient conditions of stability by means of some tangent objects in Section~\ref{Approximation}.
The results are then refined in the analytic case (Section~\ref{Analytique}) and in the linear one (Section~\ref{Lineaire}).
Section~\ref{moregeneral} is devoted to asymptotic stability for general inputs and the particular case of two globally asymptotically stable vector fields is dealt with in Section~\ref{DeuxGAS}.
Finally, in Section~\ref{Plan} a necessary and sufficient condition for a planar switched system to be GUAS is stated.
The results are illustrated by some examples in Section~\ref{Exemples}.

\section{Statement of the Problem}\label{Position}

As explained in the introduction, we deal with a finite collection $\{f_1,\ldots,f_p\}$ of smooth vector fields on $\R^d$ vanishing at the origin and assumed to share a smooth weak Lyapunov function $V: \R^d \longrightarrow \R_+$. \\
This means that $V$ is positive definite and for each $i\in\{1,\ldots,p\}$ the Lie derivative $\Li_{f_i} V$ is nonpositive, so that $V$ is nonincreasing along the solutions of $\dot{x}=f_i(x)$, $i=1,\ldots,p$. \\

An input, or switching signal, is a piecewise constant and right-continuous function $u$ from $[0,+\infty[$ into $\{1,...,p\}$. We denote by $(a_n)_{n\geq 0}$ the sequence of switching times (with $a_0=0$). In the particular case where the number of switches is finite we adopt the following convention: if $u$ is discontinuous at $a_0<\dots<a_N$ and constant from $a_N$, then we set $a_n=a_N+(n-N)$ for $n>N$. The sequence $(a_n)_{n\geq 0}$ is thus strictly increasing to $+\infty$ in any case, and $u(t)$ is constant on each interval $[a_n,a_{n+1}[$. This value will be denoted by $u_n\in \{1,...,p\}$, and the duration $a_{n+1}-a_n$ by $\delta_n$.

As the input is entirely defined by the switching times and the values taken at these instants we can write
$$
u=(a_n,u_n)_{n\geq 0}.
$$
We consider the switched system
\begin{equation}\tag{S}\label{system}
\dot{x}=f_u(x), \, x\in\R^d, u\in\{1,\dots,p\}
\end{equation}
Fixing a switching signal $u$ gives rise to the dynamical system defined in $\R^d$ by
	$$\dot{x}=f_{u(t)}(x).$$
Its solution for the initial condition $x$ will be denoted $\Phi_u(t,x)$ and $\Phi_i(t,x)$ will stand for the flow of $f_i$.\\
As we are concerned with the stability analysis at the origin of the switched system given by (\ref{system}), let us recall some stability notions :

\begin{defi}
The switched system (\ref{system}) is said to be
\begin{enumerate}
\item uniformly stable at the origin if for every $\delta>0$, there exists $\varepsilon>0$ such that for every $x\in\R^d$, $\|x\|\leq \varepsilon$, all inputs $u$ and $t\geq0$, $\|\Phi_u(t,x)\|\leq\delta$;
\item locally attractive at the origin if there exists $\delta>0$ such that for all $x\in\R^d$, $\|x\|\leq \delta$ and all inputs $u$, $\Phi_u(t,x)$ converges to the origin as $t$ goes to infinity; 
\item globally attractive at the origin if for every input $u$ and all $x\in\R^d$, $\Phi_u(t,x)$ converges to the origin as $t$ goes to infinity;
\item asymptotically stable at the origin if it is uniformly stable and locally attractive at the origin;
\item globally asymptotically stable at the origin if it is uniformly stable and globally attractive at the origin.
\item globally uniformly asymptotically stable (GUAS) at the origin if it is uniformly stable and globally uniformly attractive at the origin i.e. 
	$$\forall \varepsilon,\delta>0, \exists T>0 \quad\text{s.t.}\quad \|x\|\leq \delta \Longrightarrow \forall t\geq T \,\,  \|\Phi_u(t,x)\|\leq \varepsilon. $$
\end{enumerate}
\end{defi}

It is a well-known fact  that the set $\{f_1,\ldots,f_p\}$ being finite, hence compact, global asymptotic stability is equivalent to GUAS (see \cite{Angeli04uniformglobal}).\\   
The vector fields are not necessarily forward complete but $V$ being a  weak common Lyapunov function, the switched system is uniformly stable, and all solutions of (\ref{system}) starting in some neighborhood of the origin are forward bounded, hence defined for all positive times (however the $f_i$'s are forward complete whenever $V$ is radially unbounded, see Proposition~\ref{limitset}).\\
Before concluding this section, let us define some classes of inputs that we will consider in the sequel.

\begin{defi}\label{chaotic} (see \cite{BaldeJouan})
An input $u$ is said to be chaotic if there exists a sequence $([t_k,t_k+\tau])_{k\in\N}$ of intervals which satisfies the following conditions:
\begin{enumerate}
\item $t_k \longrightarrow_{k\to +\infty} +\infty$ and $\tau>0$;
\item For all $\varepsilon>0$ there exists $k_0\in\N$ such that for all $k\geq k_0$, the input $u$ is constant on no subinterval of $[t_k,t_k+\tau]$ of length greater than or equal to $\varepsilon$.
\end{enumerate}
An input that does not satisfy these conditions is called a nonchaotic input.
\end{defi}

In the litterature, one often encounters inputs with dwell-time (an input $u$ is said to have a dwell-time $\delta>0$ if for each $n\in\N$, $\delta_n\geq\delta$).  Such an input is clearly nonchaotic.

\begin{defi}
An input $u$ is said to satisfy the assumption $H(i)$ if there exist a subsequence $(a_{n_k})_{k\in\N}$ and $\delta>0$ such that
	$$\forall k\geq0, \ u_{n_k}=i \ \text{and} \ \delta_{n_k}\geq \delta.$$
It is said to be regular if it is nonchaotic and satisfies the assumption $H(i)$ for $i=1,\ldots,p$.
\end{defi}

\section{Limit Sets and Convergence Results}\label{Convergence}

For an input $u$ and an initial condition $x\in \R^d$, we denote by $\Omega_u(x)$ the set of $\omega$-limit points of $\{\Phi_u(t,x);\ t\geq 0\}$, that is the set of limits of sequences $(\Phi_u(t_k,x))_{k\geq 0}$, where $(t_k)_{k\geq 0}$ is increasing to $+\infty$.

Since $V$ is a Lyapunov function for each $f_i$, $i=1,\ldots,p$, it is clear that $t\in\R_+ \longmapsto V(\Phi_u(t,x))$ is nonincreasing, and we have the following result:

\begin{prop}\label{limitset}
There exists a neighborhood $\mathcal{W}$ of the origin such that for all $x\in\mathcal{W}$ and all input $u$, the solution starting at $x$ is forward bounded and the $\omega$-limit set $\Omega_u(x)$ is a compact and connected subset of a level set $\{V=R\}$ for some $R\geq0$. 

In particular if the Lyapunov function is radially unbounded ($i.e. \ V(x)\longrightarrow +\infty$ as $\|x\|$ goes to $+\infty$), we can choose $\mathcal{W}$ to be equal to $\R^d$.\
\end{prop}

\demo
$V$ has a positive minimum $m$ on the unit sphere.\\
Let us fix $r\in]0,m[$ and let $\mathcal{W}$ be the connected component of $\{V\leq r\}$ which contains $0$ (note that $\mathcal{W}$ is compact).
In the case where $V$ is radially unbounded $r$ can be chosen arbitrarily in $]0,+\infty[$.\\
Let $x\in\mathcal{W}$ and $u$ be an input.\\ 
Since $V$ is nonincreasing along the solution $\Phi_u(\cdot,x)$, the trajectory $\Phi_u(t,x)$, $t\geq0$, is included in $\mathcal{W}$, and is forward bounded. \\
It follows that $\Omega_u(x)$ is a compact and connected subset of $\mathcal{W}$ (see \cite{Marle}). 

It remains to show that $\Omega_u(x)$ is contained in a level set of $V$.\\
The map $t\in\R_+ \longmapsto V(\Phi_u(t,x))$ being nonincreasing and nonnegative, has a limit $R\geq0$ as $t$ goes to $+\infty$.\\
It is clear that $\Omega_u(x)\subseteq \{V=R\}$.
\findemo

\textit{From now on $\mathcal{W}$ will be a connected neighborhood of the origin satisfying the properties of Proposition~\ref{limitset}.\\
}

Next, we prove a uniform convergence result that will be used many times in the sequel and  which is based on the notion of weak convergence of the inputs.

\begin{prop}\label{uniform}
Let $x\in\mathcal{W}$, $u$ an input,  $(t_k)_{k\in\N}\subset\R_+$ a sequence increasing to $+\infty$ and $\tau>0$.\\
Consider the sequence $(\varphi_k)_{k\in\N}$ defined by
	$$\varphi_k(t)=\Phi_u(t_k+t,x)	\,\,, \,t\in[0,\tau], \,k\in\N.$$
Then $(\varphi_k)_{k\in\N}$ admits a subsequence that converges uniformly to an absolutely continuous function $\varphi: [0,\tau] \longrightarrow \Omega_u(x)$.\\
Moreover $\varphi$ satisfies
	$$\varphi(t)=\varphi(0)+\int_0^t \sum_{i=1}^p w_i(s)f_i(\varphi(s)) \D s , \,\, \forall t\in[0,\tau]$$	where $w_i:[0,\tau] \longrightarrow [0,1]$, $i=1,\ldots,p$, are measurable functions such that $\sum_{i=1}^p w_i=1$ almost everywhere.	
\end{prop}

\begin{rema}
Proposition~\ref{uniform} shows that a "limit trajectory" is a solution of the convexified system.
\end{rema}

\demo
The switched system $\dot{x}=f_{u(t)}(x)$ can be rewritten as the affine control one
	$$\dot{x}=\sum_{i=1}^p v_i(t)f_i(x),$$
where $v:[0,+\infty[ \longrightarrow \mathcal{B}$, and $\mathcal{B}=\{e_1,\ldots,e_p\}$ is the canonical basis of $\R^p$.\\
Consider the following functions:
	$$w^k(t)=v(t_k+t) \,\,\, \text{and} \,\, \, \varphi_k(t)=\Phi_v(t_k+t,x)$$		
for $t\in[0,\tau]$, $k\in\N$.\\
Since the solution $\Phi_u(\cdot,x)$ is forward bounded, up to a subsequence $(\varphi_k(0))_{k\in\N}$ converges to a limit $\ell\in\Omega_u(x)$ and by Proposition 10.1.5 of \cite{Sontag}, $(w^k)_{k\in\N}$ converges weakly to a measurable function $w: [0,\tau] \longrightarrow \mathcal{U}$ where $\mathcal{U}$ is the convex hull of $\mathcal{B}$ (this is equivalent to saying that $w_i\geq0$ and $\sum_{i=1}^p w_i=1$ almost everywhere). 

Since $\Phi_v(t_k+t,x)=\Phi_{w^k}(t,\varphi_k(0))$, the conclusion comes from Theorem 1 of \cite{Sontag}. 
\findemo

The following lemma will also be useful in the study of the stability for general inputs (Section~\ref{moregeneral}). The notations are the same as above.

\begin{lem}\label{faible}
Let $(w^k)_{k\in\N}\subset \Li^\infty([0,\tau],\mathcal{B})$ be a weakly convergent sequence and $w\in\Li^\infty([0,\tau],\mathcal{U})$ its limit.\\
The following statements are equivalent:
\begin{enumerate}
\item $m(\{t\in[0,\tau];\,\, w^k(t)=e_i\}) \longrightarrow 0$ (resp. $\tau$)

\item $w_i=0$ (resp. 1) almost everywhere.
\end{enumerate}
where $m$ stands for the Lebesgue's measure on the real line.
\end{lem}

\demo
It suffices to write
\begin{align}
	m(\{t\in[0,\tau];\,\, w^k(t)=e_i\})
	&=\int_0^\tau \bold{1}_{\{w^k=e_i\}}(t) \D t		\notag \\
	&=\int_0^\tau \bold{1}_{\{w_i^k=1\}}(t) \D t		\notag \\
	&=\int_0^\tau w_i^k(t) \D t		\notag
\end{align}
and to conclude using the weak convergence.
\findemo

\section{Stability Results for Regular Inputs}\label{Reguliere}

In this section we state and prove one of the main results of the paper; but first we introduce very important subsets of $\R^d$ and establish relations between them and the $\omega-$limit sets.\\
For each $i\in\{1,\ldots,p\}$ we define the set $\M_i$ by
$$
\M_i=\{x\in\R^d;\ \forall k\geq1, \,\, \Li_{f_i}^k V(x)=0\},
$$
where $\Li_{f_i}^k V= \underbrace{\Li_{f_i}\ldots\Li_{f_i}}_{k\text{ times}} V$.

The following lemma is the centerpiece of the part of the paper concerned with nonchaotic inputs because it enables us to establish two fundamental facts. On the one hand the $\omega$-limit sets intersect $\M_i$ 
for an input satisfying Assumption $H(i)$ (see Corollary~\ref{cor1}) and on the other hand, for a nonchaotic input every $\omega$-limit set is contained in the union of the $\M_i$'s (see Proposition~\ref{nonchaotic}). A fortiori the $\omega$-limit sets for regular inputs satisfy both statements. 

\begin{lem}\label{lem1}
Let $x\in\mathcal{W}$, $u$ an input and $(t_k)_{k\in\N}\subset \R_+$ a sequence increasing to $+\infty$ such that $(\Phi_u(t_k,x))_{k\in\N}$ converges to a limit $\ell\in\Omega_u(x)$.\\
If there exists $\tau>0$ such that $u$ takes the same value $i\in\{1,\ldots,p\}$ on each interval $[t_k,t_k+\tau]$
then the limit $\ell$ belongs to $\M_i$.
\end{lem}

\demo
For each $k\in\N$, $\Phi_u(t_k+\tau,x)=\Phi_i(\tau,\Phi_u(t_k,x))$, so that 
	$$\Phi_u(t_k+\tau,x) \longrightarrow_{k \to +\infty} \Phi_i(\tau,\ell).$$
It follows that $\ell$ and $\Phi_i(\tau,\ell)$ are two $\omega$-limit points of $\Omega_u(x)$, and by Proposition~\ref{limitset}, that
	$$V(\Phi_i(\tau,\ell))=V(\ell).$$
Since $V$ is nonincreasing along the trajectory $\Phi_i(\cdot,\ell)$, we get
$$\forall t\in[0,\tau], \,\, V(\Phi_i(t,\ell))=V(\ell).$$
This implies that for each $k\geq1$,
	$$\Li_{f_i}^k V(\ell)=\dfrac{d^k}{dt^k} V(\Phi_i(t,\ell))_{\big|_{t=0}}=0.$$
Finally $\ell\in\M_i$.
\findemo

\begin{cor}\label{cor1}
If an input $u$ satisfies Assumption $H(i)$, then 
	$$\forall x\in\mathcal{W}, \ \Omega_u(x)\cap\M_i\neq \emptyset.$$
\end{cor}

\begin{prop}\label{nonchaotic}
Let $x\in\mathcal{W}$ and let $u$ be a nonchaotic input. Then
	$$\Omega_u(x)\subseteq \bigcup_{i=1}^p \M_i.$$
\end{prop}

\demo
Assume that
	$$\Omega_u(x)\not\subseteq \bigcup_{i=1}^p \M_i.$$		
Let $\ell\in\Omega_u(x)\setminus \bigcup_{i=1}^p \M_i$ and let $(t_k)_{k\in\N}\subset\R_+$ be a sequence increasing to $+\infty$ such that
	$$\Phi_u(t_k,x) \longrightarrow_{k\to +\infty} \ell.$$
Let $T>0$ and set
	$$\varphi_k(t)=\Phi_u(t_k+t,x), \,\, t\in[0,T], \, k\in\N.$$	
By Proposition~\ref{uniform} and up to a subsequence, we can assume that $(\varphi_k)_{k\in\N}$ converges uniformly on $[0,T]$ to a continuous function $\varphi: [0,T] \longrightarrow \Omega_u(x)$.\\
Since $\varphi$ is continuous and $\bigcup_{i=1}^p \M_i$ is a closed subset, we can choose $T>0$ so small that
	$$\varphi([0,T])\cap \bigcup_{i=1}^p \M_i=\emptyset.$$
The input $u$ being nonchaotic, there exist  a number $\tau>0$ and a sequence $(s_k)_{k\in\N}$ such that for all $k\in\N$,
\begin{itemize}
\item $[s_k,s_k+\tau]\subseteq[t_k,t_k+T]$;
\item $u$ is constant on $[s_k,s_k+\tau]$.
\end{itemize}
Up to a subsequence, we can assume that there exists $j\in\{1,\ldots,p\}$ such that $u$ is equal to $j$ on each interval $[s_k,s_k+\tau]$, and the sequence $(\varepsilon_k)_{k\in\N}=(t_k-s_k)_{k\in\N}$ converges to a limit $\varepsilon\in[0,T]$.\\
By continuity of $\varphi$ and by uniform convergence of $(\varphi_k)_{k\in\N}$, the sequence $(\varphi_k(\varepsilon_k))_{k\in\N}$ converges to $\varphi(\varepsilon)$. According to Lemma~\ref{lem1}  the point $\varphi(\varepsilon)$ belongs to $\M_j$ , which contradicts
	$$\varphi([0,T])\cap \bigcup_{i=1}^p \M_i=\emptyset.$$
\findemo

In \cite{Bacciotti20051109} the authors showed that for all initial conditions $x$ in $\mathcal{W}$ and for every input $u$ with dwell-time, the $\omega$-limit set $\Omega_u(x)$ is contained in the union of all the compact and weakly invariant sets which are included in $\mathcal{W}\cap \bigcup_{i=1}^p \{\Li_{f_i} V=0\}$ (a compact set $M$ is said to be weakly invariant if for each $x\in M$, there exists some $i\in\{1,\ldots,p\}$ and $\delta>0$ such that for either all $t\in[0,-\delta]$ or all $t\in[0,\delta]$, $\Phi_i(t,x)\in M$).\\
Proposition~\ref{nonchaotic} is therefore a generalization of Bacciotti and Mazzi's result because the set $\bigcup_{i=1}^p \M_i$ is in general strictly contained in the set defined in \cite{Bacciotti20051109}. Moreover it applies to a larger class of inputs.

\begin{ex}\label{ridicule}
(see \cite{Bacciotti20051109,MancillaAguilar2006376})
Consider the linear switched system consisting in
	$$f_1(x)=\begin{pmatrix} -x_1-x_2 \\ x_1 \end{pmatrix}$$
and
	$$f_2(x)=\begin{pmatrix} -x_1 \\ -x_2\end{pmatrix}.$$	
$V(x)=x_1^2+x_2^2$ is a weak common Lyapunov function since
	$$\Li_{f_1} V(x)=-2x_1^2$$
and 
	$$\Li_{f_2} V(x)=-2x_1^2-2x_2^2.$$
The set $\{\Li_{f_1} V=0\}\cup\{\Li_{f_2} V=0\}=\{x_1=0\}$ is invariant for the flow of $f_2$. Bacciotti and Mazzi's result says that for inputs with dwell-time, the $\omega$-limit sets are contained in $\{x_1=0\}$.\\
 On the other hand, computing $\Li_{f_1}^3 V$ gives $\M_1\cup\M_2=\{0\}$. According to Proposition~\ref{nonchaotic}, the switched system is asymptotically stable for every nonchaotic input (recall that the class of nonchaotic inputs is larger than the class of inputs with dwell-time).\\
 Notice that the set $\{x_1=0\}$ is what we define as $\K_1 \cup \K_2$ (see Section~\ref{moregeneral}).
\end{ex}

Now we give a sufficient condition of asymptotic stability for regular inputs.

\begin{theo}\label{regstable}
If the sets $\M_i$ satisfy the condition
\begin{equation} \tag{C}
\begin{array}{l} \label{condition} 
\text{there exists a neighborhood $\V$ of the origin such that, for all $R>0$,}\\
\text{no connected component of the set $\{V=R\}\cap\mathcal{V}\cap\bigcup_{i=1}^p \M_i$}\\
\text{intersects all the $\M_i$'s,}
\end{array} 
\end{equation}
then for every regular input $u$, the switched system is asymptotically stable.\\
In particular if $\V=\mathcal{W}=\R^d$ the switched system is globally asymptotically stable for every regular switching signal.
\end{theo}

\demo
Since (\ref{system}) is uniformly stable at the origin, it suffices to check its attractivity for regular inputs.\\
Let $\V$ be a neighborhood of the origin satisfying Condition (\ref{condition}).
We can assume that $\V$ is included in $\mathcal{W}$ and is the connected component of $\{V\leq r\}$ containing $0$ for some $r>0$.

Let $x\in\V$ and let $u$ be a regular input.\\
Since $V$ is nonincreasing along the solution $\Phi_u(\cdot,x)$, we have $\Omega_u(x)\subseteq\V$.\\
We know from Proposition~\ref{limitset} that $\Omega_u(x)$ is a connected subset of a level set $\{V=R\}$ for some $R\geq 0$.
Let us show that $R=0$.\\
Since $u$ is a regular input, the two following facts hold:
\begin{enumerate}
\item  $u$ being nonchaotic, $\Omega_u(x)\subseteq \bigcup_{i=1}^p \M_i$ by Proposition~\ref{nonchaotic};
\item As $u$ satisfies $H(i)$ for $i=1,\ldots,p$, $\Omega_u(x)\cap\M_i\neq\emptyset$ for $i=1,\ldots,p$ by Corollary~\ref{cor1}.
\end{enumerate}
Therefore $\Omega_u(x)$ is a connected subset of $\{V=R\}\cap\V \cap\bigcup_{i=1}^p \M_i$ which intersects all the $\M_i$'s.
Condition (\ref{condition}) implies $R=0$.
\findemo
 
\begin{rema}
Note that Condition (\ref{condition}) cannot be fulfilled if the origin is not isolated in $\bigcap_{i=1}^p \M_i$. This provides a necessary condition for (\ref{condition}) to hold.
\end{rema}

\section{Approximation of the $\M_i$'s}\label{Approximation}

In general the sets $\M_i$ are difficult to compute and Condition (\ref{condition}) may be hard to check directly.\\
It is of interest to approximate the $\M_i$'s by some "tangent" sets easier to compute like linear subspaces or cones.

However the sets $\M_i$ are not necessarily submanifolds. Moreover $\D \Li_{f_i}^k V(0)$ may be equal to zero, it is for instance the case for $k=1$ because $\Li_{f_i} V$ has a maximum at $0$.\\
For these reasons we define $l_k$ as the smallest integer (which depends on $i$) such that $\D^{l_k} \Li_{f_i}^k V(0)\neq0$. If such an integer  does not exist we set $l_k=1$.\\
 Note that $l_1$ is necessarily even since $\Li_{f_i} V$ has a maximum at $0$.

For each $i\in\{1,\ldots,p\}$ consider the set
	$$M_i=\bigcap_{k \geq 1} \ker \D^{l_k} \Li_{f_i}^k V(0),$$
where
	$$\ker \D^{l_k} \Li_{f_i}^k V(0) = \{y\in\R^d ; \,\,\D^{l_k} \Li_{f_i}^k V(0)\cdot [y]_{l_k}=0\}$$
($[y]_{l_k}$ stands for $(\underbrace{y,\ldots,y}_{k \text{ times}})$).\\
The sets $\ker \D^{l_k} \Li_{f_i}^k V(0)$ are cones but not linear subspaces in general.
	
\begin{lem}
For each $i\in\{1,\ldots,p\}$, we have
 $$\M_i=\{x\in\R^d;\ \forall k\geq0, \,\, \Li_{f_i}^{2k+1} V(x)=0\}.$$
\end{lem}

\demo
Let $x\in\R^d$. For $t\geq0$ small enough and all $N\in\N^*$,
$$V(\Phi_i(t,x))=V(x)+\sum_{k=1}^N \Li_{f_i}^k V(x) \dfrac{t^k}{k!}+o(t^N).$$
Since $V$ is nonincreasing along solutions, the smallest integer $k$ (if it exists) such that $\Li_{f_i}^k V(x)\neq 0$ is odd, and $\Li_{f_i}^k V(x) < 0$ for this $k$.\\
It follows that 
	$$\bigcap_{k\geq0} \{\Li_{f_i}^{2k+1}V=0\}\subseteq \bigcap_{k\geq0} \{\Li_{f_i}^{2k+2}V=0\},$$
and finally that
	$$\M_i=\bigcap_{k\geq0} \{\Li_{f_i}^{2k+1}V=0\}.$$	
\findemo

\begin{lem}\label{lem2}
We have the following results:
\begin{enumerate}
\item Let $(x_n)_{n\geq 0}\subset \M_i\setminus\{0\}$ a sequence converging to the origin.\\
If $(\dfrac{x_n}{\|x_n\|})_{n\geq1}$ converges to a limit $x$, then $x\in M_i$.\\
In other words, every tangent vector to $\M_i$ at $0$ is in $M_i$.

\item If $M_i\cap M_j=\{0\}$ then $0$ is isolated in $\M_i\cap\M_j$ i.e. there exists a neighborhood $\V$ of the origin such that $\M_i\cap\M_j\cap \V=\{0\}$.

\item If $\M_i$ is a submanifold then $\mathcal{T}_0 \M_i\subseteq M_i$. In this case $M_i$ can be replaced by $\mathcal{T}_0 \M_i$.

\item $\ker \D^2 \Li_{f_i} V(0)$ is a linear subspace.
\end{enumerate}
\end{lem}

\demo
\begin{enumerate}
\item For all $k,n\geq 1$, 
	$$\Li_{f_i}^k V(x_n) = \Li_{f_i}^k V(0)+\dfrac{1}{l_k!}\D^{l_k} \Li_{f_i}^k V(0)\cdot [x_n]_{l_k}+o(\|x_n\|^{l_k}).$$
Since $x_n,0\in\M_i$, one has $\Li_{f_i}^k V(x_n)=\Li_{f_i}^k V(0) $ and
	$$\D^{l_k} \Li_{f_i}^k V(0)\cdot \bigg[\dfrac{x_n}{\|x_n\|}\bigg]_{l_k}=\dfrac{o(\|x_n\|^{l_k})}{\|x_n\|^{l_k}}.$$
Taking the limit as $n$ goes to $+\infty$, we get	
	$$\D^{l_k} \Li_{f_i}^k V(0)\cdot [x]_{l_k}=0.$$	
It follows that $x\in M_i$.		

\item If $0$ is not isolated in $\M_i\cap\M_j$, then there exists a sequence $(x_n)_{n\geq1}$ in $\M_i\cap\M_j\setminus\{0\}$ converging to $0$. \\
Up to a subsequence, $(\dfrac{x_n}{\|x_n\|})_{n\geq1}$ converges to a limit $x\neq0$ which belongs to $M_i\cap M_j$ by the previous item.

\item Since $M_i$ contains all tangent vectors to $\M_i$ at the origin, it contains the tangent space to $\M_i$ at $0$.

\item $\Li_{f_i} V$ has a maximum at $0$ so $\D^2\Li_{f_i} V(0)$ is a negative semi-definite symmetric bilinear form. It follows that the set of all isotropic vectors of $\D^2\Li_{f_i} V(0)$ is a linear subspace. 
\end{enumerate}
\findemo

\begin{theo}\label{theo2}
If no connected component of $\mathcal{S}(1)\cap \bigcup_{i=1}^p M_i$ intersects all the $M_i$'s then Condition $(C)$ is fulfilled.
\end{theo}
	
\demo
Let $I$ be the set of all indices $i$ for which there exists a continuous path in $\mathcal{S}(1)\cap \bigcup_{j=1}^p M_j$ joining $M_1$ to $M_i$, and let $J$ be the complementary set of $I$ in $\{1,\ldots,p\}$.\\
Necessarily $J$ is nonempty and for all $i\in I$ and $j\in J$, $M_i \cap M_j=\{0\}$.\\
By Lemma~\ref{lem2}, there exists a neighborhood of the origin $\V$ such that, for all $i\in I$ and $j\in J$, $\V \cap \M_i \cap \M_j=\{0\} $.\\
It is clear that Condition $(C)$ is satisfied.
\findemo

Now we give an example of a vector field for which $\M$ is hard to compute but not  $M$.
\begin{ex}
Let $\alpha: \R^2 \longrightarrow \R_-$ be a smooth nonpositive map such that $\alpha(0)=0$, which implies that $\D \alpha(0)=0$. Consider the vector field $f$ defined by $f(x)=\alpha(x)x$ and the weak Lyapunov function $V(x)=\dfrac{x_1^2+x_2^2}{2}$.\\
A straightforward computation gives $\M=\{ \alpha=0 \}$ (this set can be any closed subset of $\R^2$ by Borel's Theorem). One also has $\D^k \Li_f V(0)=0$ for $k=1,2,3$ so that $M\subseteq \ker  \D^4 \Li_f V(0)$. 
But 
\begin{align}
\dfrac{1}{12}\D^4 \Li_f V(0)\cdot [x]_4
&= \dfrac{\partial^2 \alpha}{\partial x_1^2}(0)x_1^4+\dfrac{\partial^2 \alpha}{\partial x_2^2}(0) x_2^4+2\dfrac{\partial^2 \alpha}{\partial x_1 \partial x_2}(0)(x_1^3 x_2+x_1 x_2^3)		\notag \\
&+(\dfrac{\partial^2 \alpha}{\partial x_1^2}(0)+\dfrac{\partial^2 \alpha}{\partial x_2^2}(0))x_1^2 x_2^2. \notag
\end{align}
We can choose $\alpha(x)=x_1^2-x_2^2+\beta(x)$ with $\beta(x)=o(\|x\|^2)$ too complicated to compute $\M=\{ \alpha=0 \}$. However, the previous computation gives $M\subseteq \{x_2= \pm x_1\}$. 
\end{ex}
	
\section{The Analytic Case}\label{Analytique}

In this section we assume the vector fields and the Lyapunov function to be analytic. In that case the set $\M_i$ can also be defined as
	$$\M_i=\{x\in\R^d; \,\, \forall t\in I_x ,\,\,\, V(\Phi_i(t,x))=V(x)  \},$$
where $I_x$ is the maximal interval of definition of the solution $\Phi_i(\cdot,x)$.	
It is shown in the next proposition that whenever $\V_i$ is a centre manifold for $f_i$ (see \cite{Carr} for the definition and basic properties of centre manifolds), then locally around the origin the inclusion $\M_i\subseteq \V_i$ holds and enables us to give sufficient conditions for asymptotic stability.  Thanks to the principle of approximation at any order of centre manifolds (see \cite{Carr}), these conditions are checkable in practice.

\begin{prop}\label{centre}
Let $\V_i$ be a centre manifold for $f_i$. Then there exists a neighborhood $\mathcal{O}$ of the origin such that $$\M_i\cap \mathcal{O} \subseteq \V_i.$$
\end{prop}	
	
\demo
Since $f_i$ is stable at the origin, its differential $\D f_i(0)$ at $0$ has no positive eigenvalue. Therefore we can choose a basis such that $f_i$ writes
	$$f_i(x)=Ax+g(x)$$
where $A\equiv \D f_i(0)$ is a block diagonal matrix consisting in two blocks, a $n\times n$ block $A_1$, the eigenvalues of which lies on the imaginary axis, and a $m\times m$ Hurwitz block $A_2$.\\ According to this decomposition we will write each $x$ in $\R^d$ as $x=(x_1,x_2)\in \R^n\times \R^m$.\\ 	
By definition $\V_i$ is locally the graph of a smooth map $h: \R^n \longrightarrow \R^m$ defined in a neighborhood of the origin such that $h(0)=0$ and $\D h(0)=0$ (see \cite{Carr}).

Let $\mathcal{O}$ be the connected component of $\{V\leq r\}$ containing the origin for some $r>0$.
According to Lemma 1, page 20 of \cite{Carr}, and for $r$ small enough, there exist two numbers $C,\mu>0$ such that for all $t\geq0$ and all $x=(x_1,x_2)\in \mathcal{O}$,
	$$\|x_2(t)-x_1(t)\| \leq Ce^{-\mu t} \|x_2-h(x_1)\|$$
where $x(t)=\Phi_i(t,x)$.\\
We can also assume that $\mathcal{O}$ is bounded and that $h$ is defined on $\mathcal{O}\cap (\R^n\times \{0\})$.
The map $h$ being continuous on the compact set $\mathcal{O} \cap (\R^n\times \{0\})$, there exists a number $M>0$ such that
	$$\forall x\in\mathcal{O}, \,\,  \|x_2-h(x_1)\| \leq M.$$
	
Let $x=(x_1,x_2)\in\M_i\cap\mathcal{O}$.
In order to prove that $x\in\V_i$, we show that $x_2=h(x_1)$.
By analycity $V(\Phi_i(t,x))=V(x)\leq r$ for all $t\in I_x$. It follows that the trajectory $\{\Phi_i(t,x);\,\, t\in I_x\}$ cannot leave the compact set $\mathcal{O}$ in backward time (it is therefore defined for all $t\leq 0$).\\
For each $k\in\N$, set $x^k=\Phi_i(-k,x)$ and $C_1=CM$. Then for all $t\geq0$,
\begin{align}
	\|x_2^k(t)-h(x_1^k(t))\|
	&\leq Ce^{-\mu t}\|x_2^k-h(x_1^k)\|	\notag	\\
	&\leq C_1e^{-\mu t},		\notag
\end{align} 	
 and replacing $t$ by $k$ we get
\begin{align}
	\|x_2-h(x_1)\|
	&= \|x_2^k(k)-h(x_1^k(k))\|		\notag	\\
	&\leq C_1 e^{-\mu k}	 \longrightarrow_{k\to +\infty} 0.	\notag	
\end{align}
Finally $x_2=h(x_1)$ and $x\in\V_i$.
\findemo	

From Proposition~\ref{centre}, we deduce that all tangent vectors to $\M_i$ at the origin are contained in the tangent space to $\V_i$ at the origin which will  be denoted by $V_i$. Consequently in the analytic case we get the following results.

\begin{theo}\label{tousegaux}
If one of the following conditions is satisfied:
\begin{enumerate}
\item there exists a neighborhood $\mathcal{O}$ of the origin such that for all $R>0$, no connected component of $\{V=R\}\cap\mathcal{O} \cap \bigcup_{i=1}^p \V_i$ intersects all the $\V_i$'s; 
\item no connected component of $\mathcal{S}(1)\cap \bigcup_{i=1}^p V_i$ intersects all the $V_i$'s;
\end{enumerate} 
then Condition $(C)$ is fulfilled, and the switched system is asymptotically stable for every regular input.
\end{theo}

\demo
$1. \Longrightarrow (C)$ follows from Proposition~\ref{centre}.\\
$2. \Longrightarrow 1.$ From Proposition~\ref{centre} again we get $M_i\subseteq V_i$ for $i=1,\ldots,p$. Therefore the assumptions of Theorem~\ref{theo2} are satisfied as soon as no connected component of $\mathcal{S}(1)\cap \bigcup_{i=1}^p V_i$ intersects all the $V_i$'s. 
\findemo

\begin{rema}
\begin{enumerate}
\item The linear subspaces $V_i$ are in general quite easy to compute because they are characterized by the linear part of the $f_i$'s.
\item In the case where $\V_i\subseteq \K_i$ ($\K_i$ stands for the set $\{\Li_{f_i} V=0\}$, see Section~\ref{moregeneral}) then $\M_i=\V_i$ and consequently $M_i=V_i$. Indeed, any locally positively invariant set included in $\K_i$ is contained in $\M_i$.
\end{enumerate}
\end{rema}

As well as in \cite{BaldeJouan} we deduce from Theorem~\ref{tousegaux} some geometric sufficient conditions for asymptotic stability. They deal with linear objects and are quite easily checkable.

\begin{cor} \label{prop5}
Assume that
	$$\bigcap_{i=1}^p V_i=\{0\}.$$
Then, the switched system is asymptotically stable for all regular switching signals as soon as one of the following conditions holds:
\begin{enumerate}
\item there exists $i$ such that $\dim V_i=0$;
\item there exists $i$ such that $\dim V_i=1$ and $V_i\subseteq V_j \Longrightarrow V_i=V_j$;
\item $p=2$;
\item $p>2$ and $\dim (\sum_{i=1}^p V_i)> \sum_{i=1}^p \dim V_i -p+1$. 
\end{enumerate}	
\end{cor}

\demo
Each condition implies Condition ($C$) (see \cite{BaldeJouan}).
\findemo

\begin{ex}
Let us consider the vector field $f$ defined in $\R^{3+n+m}$by
$$
\left\{
\begin{array}{ll}
\dot{x}_1 & =-x_2+y^2\varphi(z)\\
\dot{x}_2 & =x_1-y\psi(z)\\
\dot{y}   & =-y-x_1y\varphi(z)+x_2\psi(z)\\
\dot{z}   & =Bz+(x_1^2+x_2^2+y^2)Bz
\end{array}
\right.
$$
where $\varphi$ and $\psi$ are analytic functions with $\psi(0)= 0$, the variable $z$ belongs to $\R^{n+m}$, and $B$ is a Hurwitz matrix of the following form:
$$
B=\begin{pmatrix}A&-C^T\\C&D\end{pmatrix}.
$$
We assume the $n\times n$ matrix $A$ to satisfy $A^T+A=0$, the $m\times m$ matrix $D$ to satisfy $D^T+D<0$, and the pair $(C,A)$ to be observable. These conditions ensure that $B$ is Hurwitz (see \cite{BaldeJouan2}).

The positive definite function $V$ defined by $\displaystyle V(x_1,x_2,y,z)=x_1^2+x_2^2+y^2+z^Tz$ is a weak Lyapunov function for $f$. Indeed a straightforward computation gives:
$$
\Li_fV(x_1,x_2,y,z)=-2y^2+(1+x_1^2+x_2^2+y^2)z^T(B^T+B)z.
$$
As $B^T+B\leq0$, we obtain $\Li_fV\leq 0$ and:
$$
\begin{array}{ll}
\{\Li_fV=0\} & =\{y=0\}\times \ker(B^T+B)\\
           & =\{y=0\}\times\R^n.
\end{array}
$$
On the other hand the plane $P=\{y=0,\ z=0\}$ is a centre manifold for $f$. According to Proposition \ref{centre} we obtain
$$
\mathcal{M}\subseteq P\subsetneq \mathcal{K}=P\times\R^n.
$$
Thanks to LaSalle's invariance Principle and because $z(t)$ tends to zero as $t$ goes to $+\infty$ the inclusions are global.
This example can be refined in the following way. Let us first denote by $f_0$ the vector field obtained by applying to $f$ the diffeomorphism
$$
(x_1,x_2,y,z)\longmapsto (x_1,x_2,w=y+x_1^4-x_2^2,z).
$$
The $2$-dimensional manifold of $\R^{3+n+m}$ defined by $\mathcal{V}_0=\{w=x_1^4-x_2^2,\ z=0\}$ is a center manifold for $f_0$, hence $\mathcal{M}_0\subseteq \mathcal{V}_0$.

Now let $g$ be the vector field obtained by applying the same diffeomorphism to the vector field $F$ defined by:
$$
F(x_1,x_2,y,z)=(-x_1,-x_2,-y,-z) \quad \text{and let} \quad f_1=(x_1^2+x_2^2)g.
$$
The conclusion is that the switched system defined by $f_0,f_1$ is asymptotically stable for any regular input. Indeed the positive function
$$
W((x_1,x_2,w,z)=x_1^2+x_2^2+(w-x_1^4+x_2^2)^2+z^Tz
$$
is a weak Lyapunov function for $f_0,f_1$, as shown by an easy but long computation. Moreover $\mathcal{M}_1=\{x_1=x_2=0\}$ and consequently $\mathcal{M}_0\cap\mathcal{M}_1=\{0\}$.
\end{ex}

\section{The Linear Case}\label{Lineaire}

The paper \cite{BaldeJouan} deals with the linear case but consider \textbf{quadratic} Lyapunov functions only. Actually Theorems $3$ and $4$ of \cite{BaldeJouan} are direct consequences of the results of the previous sections.\\
However we will show that our results also apply to linear systems with non-quadratic, but analytic and unbounded, Lyapunov functions.\\
So let us assume the fields $f_i$ to be linear and the Lyapunov function to be analytic and radially unbounded.\\
Fix a vector field $f_i$. For every real or complex eigenvalue $a$ of $f_i$, let us denote by $L(a)=\ker (f_i-a I)^d$ the characteristic subspace of $\mathbb{C}^d$ for this eigenvalue.
 We define
$$
\begin{array}{ll}
E_a=L(a) & \mbox{if } a\in \R\\
E_a=(L(a)\oplus L(\overline{a}))\cap\R^d & \mbox{if } a\notin \R\\
\end{array}
$$
where $\overline{a}$ is the complex conjugate of $a$, and we set
$$
\begin{array}{ll}
V_i&=\sum_{\Re(a)=0}E_a\\
V_i^S&=\sum_{\Re(a)<0}E_a
\end{array}
$$
where $\Re(a)$ stands for the real part of $a$.
It is a well known fact that as soon as the real part of one of the eigenvalues of $f_i$ is positive, or in the case where the index of one of the eigenvalues whose real part vanishes is greater than one (that is the characteristic subspace is strictly larger than the eigenspace for this eigenvalue), there exists an $x$ such that $\| \Phi_i(t,x) \|$ tends to infinity as $t$ goes to $+\infty$. \\
As by hypothesis $f_i$ is stable at the origin, the space $\R^d$ is equal to $V_i\oplus V_i^S$ and the restriction of $f_i$ to $V_i$ is semi-simple (i.e. diagonalizable in $\mathbb{C}^d$).\\
 The subspaces $V_i$ and $V_i^S$ are $f_i$-invariant, so that $V_i$ is also invariant under the flow of $f_i$. It is consequently a centre manifold for $f_i$, since the restriction of $f_i$ to $V_i^S$ is Hurwitz. \\
Moreover, according to the previous decomposition of $\R^d$, we can write $f_i(x)=\begin{pmatrix}A_1 x_1\\A_2 x_2\end{pmatrix}$ where the real parts of $A_1$'s eigenvalues vanish and $A_2$ is Hurwitz. The flow of $f_i$ is then $\Phi_i(t,x)=\begin{pmatrix}e^{tA_1}x_1\\e^{tA_2}x_2\end{pmatrix}$.

\begin{prop}\label{linear}
$\M_i=\V_i=V_i$.
\end{prop}

\demo
Let $x=x_1+x_2\approx \begin{pmatrix}x_1\\x_2\end{pmatrix}\in \M_i$.\\
In order to prove that $x\in V_i$ it suffices to show that $x_2=0$.\\
Since $A_2$ is Hurwitz, there exist $C,\mu>0$ such that for all $t\in\R$ and all $y\in V_i^S$, 
$\|e^{tA_2}y\|\leq Ce^{-\mu t}\|y\|$.\\
For each $n\geq1$, set $x^n=x_1^n+x_2^n=\Phi_i(-n,x)=\begin{pmatrix} e^{-nA_1}x_1\\e^{-nA_2}x_2\end{pmatrix}$.\\
Since $x\in\M_i$ and by analycity, $V(x^n)=V(x)$ for each $n\geq 1$. As $V$ is radially unbounded $(x^n)_{n\geq1}$ is bounded, and so is $(x_2^n)_{n\geq1}$.\\
But
\begin{align}
	\|x_2\|
	&=\|\Phi_i(n,x_2^n)\|	\notag	\\
	&=\|e^{nA_2}x_2^n\|		\notag	\\
	&\leq Ce^{-\mu n}\|x_2^n\|	\notag	\\
	&\leq \tilde{C}e^{-\mu n} \longrightarrow_{n \to +\infty} 0		\notag 
\end{align}
It follows that $x_2=0$ and $x\in V_i$.\\

Conversely, since $A_1$ is semi-simple the matrices $e^{tA_1}$ are rotation ones in some Euclidean structure on $V_i$. It is a well known fact that every point of the trajectory
$t\longmapsto e^{tA_1}$
is recurrent. In particular there exists a sequence $(t_k)_{k\geq 0}$ increasing to $+\infty$ such that
$$
e^{t_kA_1}\longmapsto_{k\mapsto+\infty} I_{V_i},
$$
where $I_{V_i}$ stands for the identity matrix of $V_i$. Therefore for all $x\in V_i$
$$
V(e^{t_kA_1}x) \longmapsto_{k\mapsto+\infty} V(x),
$$
and $V(e^{tA_1}x)$ being not increasing and by analycity we obtain $x\in\M_i$.
	
\findemo

\begin{theo}
Assume that
	$$\bigcap_{i=1}^p V_i=\{0\}.$$
Then, the linear switched system is globally asymptotically stable for all regular switching signals as soon as one of the conditions (1--4) of Corollary~\ref{prop5} holds.
\end{theo}

\section{Stability for general inputs}\label{moregeneral}

In this section we deal with general inputs, more accurately the stability results of \cite{BaldeJouan} for possibly chaotic inputs are extended to nonlinear switched systems using the convergence results of Proposition~\ref{uniform} and Lemma~\ref{faible} (see Section~\ref{Convergence}).

For each $i\in\{1,\ldots,p\}$, we introduce the closed set containing the origin
	$$\K_i=\{ x\in\R^d/ \quad \Li_{f_i} V(x)=0 \},$$
and for a given input $u$, the subset $J_u$ of $\{1,\ldots,p\}$ defined by
	$$i\in J_u \Longleftrightarrow m\{t\geq0; \quad u(t)=i\}=+\infty$$
(recall that $m$ stands for the Lebesgue's measure).	

\begin{prop}\label{inclus}
For all $x\in\mathcal{W}$ and all switching signals $u$
	$$\Omega_u(x)\subseteq \bigcup_{i\in J_u} \K_i.$$
\end{prop}

\demo
Assume by contradiction that there exists $\ell\in\Omega_u(x)$ such that $\ell\not\in \bigcup_{i\in J_u} \K_i$. Let $(t_k)_{k\geq0}\subset \R_+$ be a sequence increasing to $+\infty$ and such that
	$$\ell=\lim_{k\to +\infty} \Phi_u(t_k,x).$$
Fix $T>0$ and consider the sequence $(\varphi_k)_{k\geq0}$ defined by 
	$$\forall k\in\N, \forall t\in[0,T], \quad \varphi_k(t)=\Phi_u(t_k+t,x).$$
By Proposition~\ref{uniform} and Lemma~\ref{faible}, and up to a subsequence, $(\varphi_k)_{k\geq0}$ converges uniformly on $[0,T]$ to an absolutely continuous function $\varphi: [0,T] \longrightarrow \Omega_u(x)$ satisfying for almost all $t\in[0,T]$, 
	$$ \varphi'(t)= \sum_{i\in J_u} w_i(t) f_i\circ \varphi(t)$$
where for each $i\in J_u$, $w_i:[0,T] \longrightarrow [0,1]$ is a measurable function and $\sum_{i\in J_u} w_i=1$ almost everywhere.

Since $\varphi$ is continuous, for $T$ small enough, $\varphi([0,T])\cap \bigcup_{i\in J_u} \K_i =\emptyset$.\\
Since $\Omega_u(x)$ is contained in a level set of $V$, we get for almost all $t\in[0,T]$, 
\begin{align}
	0
	&=\dfrac{\D}{\D t} V(\varphi(t))		\notag \\
	&=\D V(\varphi(t))\cdot\varphi'(t)	\notag \\
	&=\sum_{i\in J_u} w_i(t) \D V(\varphi(t))\cdot f_i(\varphi(t))		\notag \\
	&=\sum_{i\in J_u} w_i(t) \Li_{f_i} V(\varphi(t)).	\notag	
\end{align}
But for all $t\in[0,T]$, $\varphi(t)\not\in \bigcup_{i\in J_u} \K_i$ and $\Li_{f_i}V(\varphi(t))<0$ for $i\in J_u$.

It follows that for each $i\in J_u$, $w_i=0$ almost everywhere which contradicts 
	$$\sum_{i\in J_u} w_i=1 \quad \text{almost everywhere}.$$			
\findemo

\begin{prop}\label{rencontre}
For all $x\in\mathcal{W}$ and all inputs $u$: 
	$$\forall i\in J_u, \quad \Omega_u(x)\cap \K_i\neq \emptyset.$$
\end{prop}

\demo
Assume by contradiction that there exists $i\in J_u$ such that 
	$$\Omega_u(x)\cap \K_i=\emptyset.$$
Since $\Omega_u(x)$ is compact and $\K_i$ is closed, the distance between those two sets is positive, let $\varepsilon=\dfrac{d(\K_i,\Omega_u(x))}{2}>0$.\\
Fix $T>0$ such that for each $t\geq T$, 
	$$d(\Phi_u(t,x),\Omega_u(x))\leq\varepsilon.$$	
Since $K:=\{ y\in\R^d/ \,\, d(y,\Omega_u(x))\leq \varepsilon \}$ is a compact set which does not intersect $\K_i$ there exists $\alpha>0$ such that 
	$$\forall y\in K, \quad \Li_{f_i}V(y)\leq -\alpha.$$	
For each $t\geq0$,
\begin{align}
	V(\Phi_u(t,x))
	&=V(x)+\int_0^t  \dfrac{\D}{\D s} V(\Phi_u(s,x)) \D s		\notag \\
	&=V(x)+\int_0^t  \D V(\Phi_u(s,x))\cdot f_{u(s)}(\Phi_u(s,x)) \D s.			\notag
\end{align}
As $\Omega_u(x)$ is contained in a level set $\{V=R\}$ for some $R>0$, this integral converges as $t$ goes to infinity and we get
\begin{align}
	R
	&=V(x)+\int_0^{+\infty}  \D V(\Phi_u(s,x))\cdot f_{u(s)}(\Phi_u(s,x)) \D s			\notag \\
	&\leq V(x)+\int_T^{+\infty}  \D V(\Phi_u(s,x))\cdot f_{u(s)}(\Phi_u(s,x)) \D s			\notag \\
	&\leq V(x)+\int_{\{t\geq T; \,\, u(t)=i\}}  \D V(\Phi_u(s,x))\cdot f_i(\Phi_u(s,x)) \D s		\notag \\
	&\leq V(x)-\alpha m\{t\geq T; \,\, u(t)=i\}=-\infty, 	\notag
\end{align}	
a contradiction.
\findemo

From Propositions~\ref{inclus} and~\ref{rencontre}, we deduce for general inputs a theorem similar to Theorem~\ref{regstable} for regular ones.

\begin{theo}\label{chaostable}
Let $J\subseteq \{1,\ldots,p\}$ be a nonempty set.\\
If the sets $\K_i$ satisfy the condition
\begin{equation} \tag{K}
\begin{array}{l} \label{conditionK} 
\text{there exists a neighborhood $\V$ of the origin such that, for all $R>0$,}\\
\text{no connected component of the set $\{V=R\}\cap\mathcal{V}\cap\bigcup_{i\in J} \K_i$}\\
\text{intersects all the $\K_i$, for $i\in J$}
\end{array} 
\end{equation}
then for every input $u$ such that $J_u=J$, the switched system is asymptotically stable.\\
In particular if $\V=\mathcal{W}=\R^d$ then the switched system is globally asymptotically stable for every input $u$ such that $J_u=J$.
\end{theo}

Using the same technique of linearization as in the regular input case we can give sufficient conditions for (\ref{conditionK}) to hold. Set $K_i=\ker \D^2 \Li_{f_i} V(0)$, $i=1,\ldots,p$, and recall that these sets are linear subspaces (see Lemma~\ref{lem2}).

\begin{cor} \label{LinearisationK}
Let $J\subseteq \{1,\ldots,p\}$ be a nonempty set and assume that
	$$\bigcap_{i\in J} K_i=\{0\}.$$
Then, the switched system is asymptotically stable for all  switching signals $u$ such that $J_u=J$ as soon as one of the following conditions holds:
\begin{enumerate}
\item there exists $i\in J$ such that $\dim K_i=0$;
\item there exists $i\in J$ such that $\dim K_i=1$ and $K_i\subseteq K_j \Longrightarrow K_i=K_j$ for each $j\in J$;
\item $\#J=2$;
\item $\#J>2$ and $\dim (\sum_{i\in J} K_i)> \sum_{i\in J} \dim K_i -\#J+1$. 
\end{enumerate}	
\end{cor}

\section{The case of two GAS vector fields}\label{DeuxGAS}

In this section, we are interested in asymptotic stability for all inputs of a system consisting in a pair of vector fields. They are of course assumed to be asymptotically stable at the origin (consider constant inputs) and denoted by $f_0$ and $f_1$ so that the switched system (\ref{system}) writes
	$$\dot{x}=(1-u)f_0(x)+uf_1(x), \,\, x\in\R^d, u\in\{0,1\}.$$
Since the vector fields are asymptotically stable, we can choose a neighborhood $\mathcal{W}$ of the origin such that $\forall x\in\mathcal{W}$ and for $i=0,1$, $\Phi_i(t,x)$ goes to $0$ as $t$ goes to $+\infty$. \\
We consider the set 
	$$\K=\K_0\cap\K_1$$
and we get:

\begin{prop}\label{contenuK}
For all $x\in\mathcal{W}$ and all inputs $u$, 
	$$\Omega_u(x)\cap \K\neq \emptyset.$$
Moreover, if the vector fields and the Lyapunov function are analytic then
	$$\Omega_u(x)\subseteq \K.$$
	
If the vector fields are globally asymptotically stable and $V$ is unbounded, one can take $\mathcal{W}=\R^d$.	
\end{prop}

\demo
Assume by contradiction that $\Omega_u(x)\setminus\K\neq \emptyset$.\\
Let $\ell\in\Omega_u(x)\setminus\K$ and let $(t_k)_{k\in\N}\subset \R_+$ be a sequence increasing to $+\infty$ such that
	$$\Phi_u(t_k,x) \longrightarrow_{k\to +\infty} \ell.$$
Since $\Omega_u(x)\subseteq \K_0\cup\K_1$ we can assume without loss of generality that $\ell\in\K_0\setminus\K_1$.	\\
Fix $\tau>0$ and consider the sequence $(\varphi_k)_{k\in\N}$ defined by
	$$\varphi_k(t)=\Phi_u(t_k+t,x), \,\, t\in[0,\tau].$$
By Proposition~\ref{uniform}, up to a subsequence, $(\varphi_k)_{k\in\N}$ uniformly converges on $[0,\tau]$ to an absolutely continuous function $\varphi: [0,\tau] \longrightarrow \Omega_u(x)$	satisfying almost everywhere
	$$\dot{\varphi}=(1-\lambda)f_0(\varphi)+\lambda f_1(\varphi)$$	
where $\lambda: [0,\tau] \longrightarrow [0,1]$ is a measurable function.\\
We can choose $\tau$ so small that 
	$$\varphi([0,\tau])\subseteq \K_0\setminus\K_1.$$
As $\Omega_u(x)$ is contained in a level set $\{V=r\}$ for some $r>0$, we get for almost all $t\in[0,\tau]$
\begin{align}
0
&=\dfrac{\D}{\D t} V(\varphi(t))  \notag \\
&=(1-\lambda(t))\Li_{f_0} V(\varphi(t))+\lambda(t)\Li_{f_1}(\varphi(t)).  \notag 
\end{align}
But $\varphi([0,\tau])\cap \K_1=\emptyset$, and $\Li_{f_1} V(\varphi(t))<0$ for all $t\in[0,\tau]$. \\
Consequently $\lambda(t)=0$ for almost all $t\in[0,\tau]$ and
	$$\dot{\varphi}=f_0(\varphi) \quad \text{i.e.} \quad \forall t\in[0,\tau], \,\, \varphi(t)=\Phi_0(t,l).$$
In the analytic case this implies $\forall t\geq 0$ $V(\Phi_0(t,\ell))=r$, and $\Phi_0(t,\ell)$ cannot converge to $0$ as $t$ goes to $+\infty$.

In the smooth case let us assume that $\Omega_u(x)\cap \mathcal{K}=\emptyset$. Let us first notice that the limit trajectory can be extended to any interval $[0,T[$ hence to $[0,+\infty[$. The assumption $\Omega_u(x)\cap \mathcal{K}=\emptyset$ implies then that $\varphi(t)\in \mathcal{K}_0\setminus\mathcal{K}_1$ for $t\geq 0$, hence that $\varphi(t)=\Phi_0(t,\ell)$ for $t\geq 0$. The conclusion comes from the same contradiction that $\Phi_0(t,\ell)$ cannot converge to $0$ as $t$ goes to $+\infty$.

\findemo

\begin{cor}\label{Kzero}
Assume that the vector fields $f_0$ and $f_1$ are globally asymptotically stable and that $V$ is radially unbounded.
If $\K=\{0\}$ then the switched system is $GUAS$.
\end{cor}

\section{Planar analytic switched systems}\label{Plan}

In this section we give a necessary and sufficient condition for a planar analytic switched system consisting in two globally asymptotically stable subsystems to be GUAS under our usual assumption of existence of a weak common Lyapunov function.

 Consider the analytic planar switched system
\begin{equation}\label{planarsystem}\tag{P}
\dot{x}=(1-u)f_0(x)+uf_1(x), \,\, x\in\R^2, u\in\{0,1\}
\end{equation}		
where  $f_0$, $f_1$ are two globally asymptotically stable analytic vector fields and $V$ is a weak analytic Lyapunov function.
		 
We  assume $V$ to be radially unbounded and $\D V(x)$ to vanish at $x=0$ only, which implies that for each $R>0$, the set $\{V=R\}$ is a one-dimensional submanifold.\\
It is a well-known fact that (\ref{planarsystem}) is GUAS if and only if the convexified planar switched system, obtained by replacing $u\in \{0,1\}$ by $u\in [0,1]$, is GUAS (see \cite{Angeli04uniformglobal,BoscainCharlotSigalotti}).	\\	 
In \cite{BoscainCharlotSigalotti} the authors consider the set
	$$\mathcal{Z}=\{x\in\R^2: \det(f_0(x),f_1(x))=0\}$$
where the vector fields are colinear.
Under generic conditions they give a necessary and sufficient condition for a planar nonlinear switched system to be GUAS, a necessary condition being that the set $\{x\in\mathcal{Z}, \langle f_0(x),f_1(x) \rangle<0\}$ is empty.\\
This condition turns out to be also sufficient under our hypothesis.

\begin{theo}\label{colinear}
Under our hypothesis the following statements are equivalent:
\begin{enumerate}
\item the switched system (\ref{planarsystem}) is GUAS;
\item the set $\{x\in\mathcal{Z}: \langle f_0(x),f_1(x)\rangle<0\}$ is empty;
\item the set $\{x\in\K: \langle f_0(x),f_1(x)\rangle<0\}$ is empty.
\end{enumerate}
\end{theo}

\demo
$1. \Longrightarrow 2.$ Let $y\in\{x\in\mathcal{Z}: \langle f_0(x),f_1(x)\rangle<0\} \neq \emptyset$.\\
There exists $u_0\in[0,1]$ such that $(1-u_0)f_0(y)+u_0 f_1(y)=0$. The point $y$ is therefore an equilibrium of the convexified switched system for the input $u\equiv u_0$.\\
The convexified system is not GUAS, and neither is System (\ref{planarsystem}).\\

$2. \Longrightarrow 1.$ Assume that $\{x\in\mathcal{Z}: \langle f_0(x),f_1(x)\rangle<0\}=\emptyset$.\\
Let $x\in\R^2$ and $u$ an input. Assume by contradiction that $\Omega_u(x)\neq\{0\}$.\\ 
Let $\varphi:[0,+\infty[ \longrightarrow \Omega_u(x)$ be a limit trajectory as in Proposition~\ref{uniform}.\\	
There exists a measurable function $\lambda:[0,+\infty[ \longrightarrow [0,1]$ such that
	$$\dot{\varphi}=(1-\lambda)f_0(\varphi)+\lambda f_1(\varphi) \quad \text{a.e.}$$
The subsystems and the Lyapunov function being analytic, we know from Proposition~\ref{contenuK} that $\Omega_u(x)$ is contained in $\K$ so that $\Omega_u(x)\subseteq \mathcal{Z}$. 
Indeed, for each $y\in\K\setminus \{0\}$, $\Li_{f_0} V(y)=\Li_{f_1} V(y)=0$ and this is equivalent to $f_0(y),f_1(y)\in \ker \D V(y)$. But $\D V(y)\neq 0$;  it follows that $f_0(y)$ and $f_1(y)$ are colinear and that $\K\subseteq \mathcal{Z}$. In particular, $\varphi([0,+\infty[)\subseteq \mathcal{Z}$.\\
Since $\{x\in\mathcal{Z}: \langle f_0(x),f_1(x)\rangle<0\}=\emptyset$, for all $y\in\mathcal{Z}\setminus\{0\}$
	$$f_1(y)=\alpha(y)f_0(y), \,\,\,\text{where} \,\,\, \alpha(y)=\dfrac{\|f_1(y)\|}{\|f_0(y)\|}>0.$$ 
It follows that
	$$\dot{\varphi}=(1-\lambda+\lambda\alpha(\varphi))f_0(\varphi)$$
with $1-\lambda+\lambda\alpha(\varphi)>0$, so the trajectories $\{\varphi(t): t\geq0\}$ and $\{\Phi_0(t,\varphi(0)): t \geq0\}$ are the same which contradicts the global asymptotic stability of $f_0$. Consequently the switched system is GUAS.\\

$2. \Longrightarrow 3.$ comes from the fact that $\K\subseteq \mathcal{Z}$.\\

$3. \Longrightarrow 2.$ It suffices to show that $\{x\in\mathcal{Z}: \langle f_0(x),f_1(x) \rangle<0\} \subseteq \K$.\\
Let $y\in\{x\in\mathcal{Z}: \langle f_0(x),f_1(x)\rangle<0\}$ and let $\lambda\in]0,1[$ such that $(1-\lambda)f_0(y)+\lambda f_1(y)=0$. Then $0\leq -\lambda \Li_{f_1}(y)=(1-\lambda)\Li_{f_0}(y)\leq0$, and $\Li_{f_0}V(y)=\Li_{f_1}V(y)=0$, i.e. $y\in\K$.
\findemo

\begin{cor}
If for all $x\in\R^2$ $\langle f_0(x),f_1(x)\rangle \geq0$, then (\ref{planarsystem}) is GUAS.
\end{cor} 

As the next example shows, the previous theorem is no longer true without the analycity hypothesis.

\begin{ex}
Consider the smooth vector field
	$$g(x)=\begin{pmatrix} -x_2 \\ x_1\end{pmatrix}$$
whose trajectories are circles around the origin, and the closed sets $F_1=\{ x_1x_2\geq0\}$ and $F_2=\{x_1x_2\leq0\}$.\\
By Borel's Theorem, there exists two smooth functions $\varphi_i: \R^2 \longrightarrow \R_+$, $i=1,2$, such that $F_i=\varphi_i^{-1}(\{0\})$.\\
Then we consider the vector fields
	$$f_i(x)=-\varphi_i(x)x+g(x), \,\,i=1,2$$
and the weak Lyapunov function $V(x)=x_1^2+x_2^2$.\\
It is easy to see that $\K_i=F_i$, $\K=\mathcal{Z}=\{x_1=0\}\cup\{x_2=0\}$ and for each $x\in\mathcal{Z}$, $\langle f_1(x),f_2(x)\rangle \geq0$.\\
Thus all hypothesis of Theorem~\ref{colinear} are satisfied except analycity of the vector fields. However the switched system is not GUAS since it admits periodic trajectories. 
\end{ex}

\section{Examples}\label{Exemples}

\begin{ex}Consider the following analytic vector fields on $\R^2$,
	$$f_1(x)=(x_1^2-x_2^3)^2\begin{pmatrix} -x_1 \\ -x_2\end{pmatrix}$$
and	
	$$f_2(x)=(x_2^2-x_1^3)^2\begin{pmatrix} -x_1 \\ -x_2\end{pmatrix}.$$
Then $V(x)=x_1^2+x_2^2$ is a weak common Lyapunov function which is radially unbounded,
	$$\Li_{f_1} V(x)=-2(x_1^2+x_2^2)(x_1^2-x_2^3)^2,$$ 
and	
	$$\Li_{f_2} V(x)=-2(x_1^2+x_2^2)(x_2^2-x_1^3)^2.$$
It follows that $\Li_{f_1} V(x)=0$ if and only if $x_1^2-x_2^3=0$. 
But all these points are equilibrium so that 
	$$\M_1 =\K_1= \{x_1^2-x_2^3=0\}.$$
Similarly $\Li_{f_2} V(x)=0$ if and only if $x_2^2-x_1^3=0$, and 
	$$\M_2 =\K_2= \{x_2^2-x_1^3=0\}.$$
For $0<R<2$, no connected component of $\{V=R\}\cap (\M_1 \cup \M_2)$ intersects $\M_1$ and $\M_2$. By Theorem~\ref{regstable} the switched system is asymptotically stable for all regular input.\\
Notice that we cannot conclude using centre manifolds because the linear parts of $f_1$ and $f_2$  at the origin vanish.	
\end{ex}

\begin{ex}Consider the following vector fields on $\R^2$,
	$$f_1(x)=\begin{pmatrix} 0 \\ -x_2\end{pmatrix}$$
and	
	$$f_2(x)=\begin{pmatrix} x_1x_2 \\ -x_2-x_1^2\end{pmatrix}.$$
The weak Lyapunov function $V(x)=x_1^2+x_2^2$ is radially unbounded, and
	$$\Li_{f_1} V(x)=\Li_{f_2} V(x)=-2x_2^2.$$ 
From
         $$\Li_{f_2}^3 V(x)=-4(x_1^3+x_2x_1^2+2x_2^2),$$
follows that $\M_1=\{x_2=0\}$ and $\M_2=\{0\}$, and that the switched system is globally asymptotically stable for all regular switchings.	\\
We could conclude using centre manifolds but the conclusion would be only local.\\
Indeed $\V_1=\{x_2=0\}$ and $\V_2$ is the graph of a map $x_2=-x_1^2+O(x_1^4)$ (see \cite{Carr}).\\
It follows that Condition 1 of Theorem~\ref{tousegaux} is satisfied and that the switched system is asymptotically stable for all regular inputs. However we cannot say anything about the region of attraction.         	
\end{ex}

\begin{ex}(see \cite{Aleksandrov2012127})
Now we consider
	$$f_1(x)=\begin{pmatrix} x_2^3 \\ -x_1^3-2x_2^3\end{pmatrix}$$ 
and 
	$$f_2(x)=\begin{pmatrix} -2x_1^3-x_2^3 \\ x_1^3\end{pmatrix}.$$\\
Here $V(x)=x_1^4+x_2^4$ is a weak common Lyapunov function.
Indeed,
	$$\Li_{f_1} V(x)=-8x_2^6$$
and
	$$\Li_{f_2} V(x)=-8x_1^6.$$	
It follows that $\K_1=\{x_2=0\}$ and $\K_2=\{x_1=0\}$ so that $\K=\{0\}$. 	\\
According to Corollary~\ref{Kzero}, the switched system is GUAS.
\end{ex}

\begin{ex}
(see \cite{Bacciotti20051109})
Consider the switched system consisting in the two following vector fields
	$$f_1(x)=\begin{pmatrix} -x_2 \\ x_1-x_2^k\end{pmatrix}$$
and
	$$f_2(x)=\begin{pmatrix} x_2 \\ -x_1-x_2^k\end{pmatrix}$$	
where $k$ is an odd integer.\\
The map $V(x)=x_1^2+x_2^2$ is a weak common Lyapunov function. Indeed,
	$$\Li_{f_1} V(x)=\Li_{f_2} V(x)=-2x_2^{k+1}.$$
Then $\M_1$ and $\M_2$ are contained in the set $\K_1=\K_2=\{x_2=0\}$.\\
We know from Proposition~\ref{inclus} that for all inputs $u$,
	$$\forall x\in\R^2, \,\,\Omega_u(x)\subseteq \{x_2=0\}.$$
But each limit set is connected and contained in a level set of $V$, hence reduced to a point.\\
If  $u$ satisfies Assumption $H(i)$ for at least one $i\in\{1,2\}$ then for each $x\in\R^2$, $\Omega_u(x)$ consists in an equilibrium point for the vector field $f_i$ which implies that $\Omega_u(x)=\{0\}$.\\
So the switched system is globally asymptotically stable for all inputs satisfying $H(i)$ for at least one $i\in\{1,2\}$.
\end{ex}

\begin{ex}
Consider the switched system consisting in the two following analytic vector fields
	$$f_1(x)=\begin{pmatrix} -x_1^3-x_2 \\ x_1^3\end{pmatrix}$$
and
	$$f_2(x)=\begin{pmatrix} -2x_1^3-x_2 \\ x_1^3\end{pmatrix}.$$	
The function $V(x)=x_1^4+2x_2^2$ is a weak analytic common Lyapunov function. Indeed,
	$$\Li_{f_1} V(x)=-4x_1^6$$
and
	$$\Li_{f_2} V(x)=-8x_1^6.$$
It follows that $\K=\{x_1=0\}$, so that Condition (\ref{conditionK}) is not fulfilled.\\
However the LaSalle invariance principle tells us that the vector fields are globally asymptotically stable, and on the other hand, for each $x\in\K$, $\langle f_1(x),f_2(x)\rangle =x_2^2\geq0$. According to Theorem~\ref{colinear} the switched system is GUAS.	
\end{ex}
 
\begin{ex}
Consider the three following analytic vector fields:
	$$f_1(x)=\begin{pmatrix} -2x_2+x_1x_3 \\ 
						x_1 \\
						-x_3-x_1^2               \end{pmatrix}$$

	$$f_2(x)=\begin{pmatrix} -x_1 \\
	                                               -2x_2 \\
	                                               -x_1^2 x_3^3\end{pmatrix}$$
	
	$$f_3(x)=\begin{pmatrix} -x_1+2x_3-2x_2 \cos x_1 \\
					         x_1 \cos x_1 \\
					         -x_3\end{pmatrix}$$	
and the analytic Lyapunov function $V(x)=x_1^2+2x_2^2+x_3^2$.	
It is easy to see that $V_1=\{x_3=0\}$, $V_2=\{x_1=x_2=0\}$ and $V_3=\{x_1=x_3=0\}$ so that the switched system is asymptotically stable for regular inputs by Theorem~\ref{tousegaux} (or Condition 2 of Corollary~\ref{prop5}).
We have
	$$\Li_{f_1} V(x)=-2x_3^2$$
	$$\Li_{f_2} V(x)=-2x_1^2-8x_2^2-2x_1^2 x_3^4$$
	$$\Li_{f_3} V(x)=-(x_1-x_3)^2$$
thus $\K_1=V_1$, $\K_2=V_2$ and $\K_3=\{x_1=x_3\}$.	It follows by Theorem~\ref{chaostable} that the switched system is globally asymptotically stable for all inputs $u$ such that $J_u=\{1,2,3\}$.
\end{ex} 

\newpage
\nocite*
\bibliography{biblio}
\bibliographystyle{plain}
\end{document}